\documentclass[12pt]{article}

\usepackage{amsmath,amssymb,amsthm,mathrsfs}
\usepackage{bm}
\usepackage{graphicx}
\usepackage[
  a4paper,
  left=25mm,
  right=25mm,
  top=38.5mm,
  bottom=38.5mm
]{geometry}
\usepackage{hyperref}
\hypersetup{
  colorlinks=true,
  linkcolor=blue,
  citecolor=red,
  urlcolor=blue
}
\usepackage{enumitem}

\makeatletter
\@addtoreset{equation}{section}
\@addtoreset{equation}{subsection}
\renewcommand{\theequation}{%
  \ifnum\value{subsection}=0\relax%
    \thesection.\arabic{equation}%
  \else%
    \thesubsection.\arabic{equation}%
  \fi}

\makeatother

\newtheorem{theorem}{Theorem}[section]
\newtheorem{lemma}[theorem]{Lemma}
\newtheorem{proposition}[theorem]{Proposition}
\newtheorem{corollary}[theorem]{Corollary}

\theoremstyle{definition}
\newtheorem{definition}[theorem]{Definition}
\newtheorem{remark}[theorem]{Remark}

\title{Asymptotic behavior of solutions to the one-dimensional $p$-system\\
with time-dependent non-smooth damping}

\author{Toi Fujimura\thanks{
Department of Mathematics, Faculty of Science and Technology, Tokyo University of Science,
2641 Yamazaki, Noda, Chiba 278-8510, Japan.
E-mail: \texttt{toui.fujimura@gmail.com}} \and
Motohiro Sobajima\thanks{
Department of Mathematics, Faculty of Science and Technology, Tokyo University of Science,
2641 Yamazaki, Noda, Chiba 278-8510, Japan.
E-mail: \texttt{msobajima1984@gmail.com}}}

\date{}

\begin{document}

\maketitle

\begin{abstract}
This paper is concerned with the asymptotic behavior of solutions to the Cauchy problem for the one-dimensional $p$-system with time-dependent damping.
When the damping coefficient $b(t)$ behaves like $(1+t)^{-\lambda}$ with $0\le\lambda<1$, the solutions are expected to converge to a diffusion wave.
For $b(t)=(1+t)^{-\lambda}$, previous studies have established convergence to a diffusion wave together with convergence rates using the energy method, provided that the initial perturbation and the wave strength are sufficiently small.
In this paper, we introduce two auxiliary functions and reformulate our problem into two perturbation problems adapted to the structure of the system.
This enables us to prove a result analogous to those obtained in previous studies, without assuming the differentiability of $b(t)$.
The auxiliary functions are given by the differential equations $1+g'(t)=b(t)g(t)$ and $1-h'(t)=b(t)h(t)$, respectively.
Our result shows that it suffices to assume $b\in C([0,\infty))$, $b(t)\asymp(1+t)^{-\lambda}$, and $|g'(t)|\lesssim(1+t)^{-(1-\lambda)}$ to obtain convergence to a diffusion wave.
\end{abstract}



\section{Introduction}\label{sec:introduction}

We consider the following Cauchy problem for the one-dimensional $p$-system with time-dependent damping:
\begin{equation}\label{eq:p-system}
\left\{
\begin{aligned}
&v_t-u_x=0, \\
&u_t+(p(v))_x=-b(t)u, &\qquad& (x,t)\in \mathbb{R}\times (0,\infty), \\
&(v,u)(x,0)=(v_0,u_0)(x), && x\in \mathbb{R}.
\end{aligned}
\right.
\end{equation}
Here $p\in C^\infty((0,\infty))$ satisfies $p>0$ and $p'<0$, and $b(t)>0$ describes the effect of dissipation.
In particular, we consider non-smooth $b(t)$.
Physically, this system models one-dimensional compressible flow through porous media, where $v=v(x,t)>0$ represents the specific volume, $u=u(x,t)$ the velocity, and $p=p(v)$ the pressure; a typical example is $p(v)=v^{-\gamma}$ ($\gamma \ge 1$).

We assume the far-field condition:
\begin{equation}\label{eq:far-field}
(v_0,u_0)(x)\to (v_\pm,u_\pm) \qquad \text{as} \quad x\to\pm\infty,
\end{equation}
where $v_\pm>0$ and $u_\pm$ are constants.
In this paper, we consider the case where
\begin{equation}\label{eq:end-states}
v_+ \neq v_-.
\end{equation}
The case $v_+=v_-$ has also been studied in the literature (e.g.,
\cite{Hsiao-Liu} and \cite{Nishihara} for constant damping,
\cite{Cui-Hou} for time-dependent damping,
\cite{Matsumura-Nishihara} for space-dependent damping, and
\cite{Cai} for spatiotemporal damping).
The case $v_+\ne v_-$ with spatiotemporal damping has also been studied in \cite{Cai2}.

Hereafter, we consider the case $v_+ \neq v_-$. For the constant damping case,
namely, $b(t)\equiv\alpha>0$, Hsiao--Liu~\cite{Hsiao-Liu} studied small
perturbations around a nonlinear diffusion wave using the energy method.
They showed that, provided that the initial perturbation and
the wave strength $\delta:=|v_+-v_-|+|u_+-u_-|$ are sufficiently small, the solution of
\eqref{eq:p-system}--\eqref{eq:end-states} behaves asymptotically like a
solution $(\bar v,\bar u)(x,t)$ of the following system based on Darcy's law:
\begin{equation}\label{eq:constant-profile}
\left\{
\begin{aligned}
&\bar v_t - \bar u_x = 0,\\
&(p(\bar v))_x =-\alpha \bar u.
\end{aligned}
\right.
\end{equation}
Moreover, the convergence rate $O\bigl((1+t)^{-\frac12}\bigr)$ was obtained
in both the $L^2(\mathbb{R})$ and $L^\infty(\mathbb{R})$ norms for the solutions toward the
corresponding profile constructed from $(\bar v,\bar u)$.

Under the same smallness assumptions, Nishihara~\cite{Nishihara} subsequently improved the higher-order $L^2$
decay estimates and derived the following sharper $L^\infty$ rates from
them:
\[
\|(v-\bar v)(\cdot,t)\|_{L^\infty} = O\bigl((1+t)^{-\frac34}\bigr), \qquad
\|(u-\bar u)(\cdot,t)\|_{L^\infty} = O\bigl((1+t)^{-\frac54}\bigr)
\]
as $t\to\infty$.

Here we give the observations.  For the constant damping case
$b(t)\equiv \alpha$, eliminating $u$ from \eqref{eq:p-system} yields the
nonlinear damped wave equation
\begin{equation}\label{eq:nonlinear-wave}
v_{tt}+(p'(v)v_x)_x+\alpha v_t=0, \qquad p'(v)<0.
\end{equation}
On the other hand, eliminating $\bar u$ from \eqref{eq:constant-profile}, we obtain the corresponding nonlinear diffusion equation
\[
\bar v_t+\frac{1}{\alpha}(p'(\bar v)\bar v_x)_x=0.
\]
Thus, the above results in \cite{Hsiao-Liu,Nishihara} show the diffusion
phenomenon, in which solutions to the damped wave equation behave
asymptotically like solutions to the corresponding nonlinear diffusion equation.

We now introduce a quantity $m$ which plays a crucial role in analyzing the asymptotic behavior.
Put $m=v_t+\alpha v$.
Then equation \eqref{eq:nonlinear-wave} can be written in the divergence form as
\begin{equation}\label{eq:div-form}
m_t+(p'(v)v_x)_x=0.
\end{equation}
Since $v_{tt}$ is expected to decay faster than $v_t$ as $t\to\infty$, we formally have $m_t\approx\alpha v_t$.
Thus, reformulating the system in divergence form also allows us to predict the asymptotic profile.

Next, we turn to the time-dependent case.
According to Wirth~\cite{Wirth}, for the linear damped wave equation
\begin{equation}\label{eq:linear-wave}
v_{tt}-v_{xx}+b(t)v_t=0,
\end{equation}
the diffusion phenomenon occurs in the effective case
$b(t)=(1+t)^{-\lambda}$ with $-1\le\lambda<1$.

Our quasilinear problem \eqref{eq:p-system} has been studied under the condition
\[
b(t)=(1+t)^{-\lambda}, \qquad 0\le\lambda<1.
\]
In this case,  under the same smallness assumptions, Li--Li--Mei--Zhang~\cite{Li} and Cui--Yin--Zhang--Zhu~\cite{Cui} proved convergence to a solution of the following system obtained by replacing $\alpha$ in \eqref{eq:constant-profile} with $b(t)$:
\begin{equation}\label{eq:variable-profile}
\left\{
\begin{aligned}
&\bar v_t - \bar u_x = 0,\\
&(p(\bar v))_x =-b(t) \bar u.
\end{aligned}
\right.
\end{equation}
In \cite{Li}, an argument analogous to that in \cite{Nishihara} was carried out and
$L^2$ decay estimates were obtained, from which $L^\infty$ convergence
rates toward the asymptotic profile given by \eqref{eq:variable-profile}
were derived.
In \cite{Cui}, both the $L^2$ and $L^\infty$ decay rates were further
sharpened; in particular, the following estimates were obtained:
\begin{subequations}\label{eq:cui-rates}
\begin{equation}
\|(v-\bar v)(\cdot,t)\|_{L^{\infty}}=
\left\{
\begin{aligned}
&O\bigl((1+t)^{-\frac{3(1+\lambda)}{4}}\bigr), &\qquad& 0\le \lambda<1/7, \\
&O\bigl((1+t)^{-\frac67+\eta}\bigr), && \lambda=1/7, \\
&O\bigl((1+t)^{-(1-\lambda)}\bigr), && 1/7<\lambda <1,
\end{aligned}
\right.
\end{equation}
\begin{equation}
\|(u-\bar u)(\cdot,t)\|_{L^{\infty}}=
\left\{
\begin{aligned}
&O\bigl((1+t)^{-\frac{5+\lambda}{4}}\bigr), &\qquad& 0\le \lambda<1/7, \\
&O\bigl((1+t)^{-\frac97+\eta}\bigr), && \lambda=1/7, \\
&O\bigl((1+t)^{-\frac{3(1-\lambda)}{2}}\bigr), &&1/7<\lambda<1
\end{aligned}
\right.
\end{equation}
\end{subequations}
as $t\to\infty$, where $\eta>0$ is arbitrary.
Here the critical case $\lambda=1/7$ arises, and the polynomial loss $(1+t)^\eta$ remains. 
For $1/7<\lambda<1$, Geng--Huang--Jin--Wu~\cite{GengHuangJinWu} further constructed
a time-asymptotic expansion around the diffusion wave and showed that it
provides the best asymptotic profile.

Note that differentiability of the damping coefficient $b(t)$ is required in their analysis.
Indeed, if we put $m=v_t+b(t)v$ and derive an equation similar in form to \eqref{eq:div-form}, a derivative of $b(t)$ appears as a remainder term:
\[
(v_t+b(t)v)_t+(p'(v)v_x)_x=b'(t)v.
\]

On the other hand, in Sobajima~\cite{Sobajima'23, Sobajima'25}, an energy method for the linear problem \eqref{eq:linear-wave} with $b \notin C^1$ is introduced.
The key idea is to use an auxiliary function $g(t)$ appropriately, as explained below.
Motivated by this idea, we introduce $m=m_g$ adapted to the time-dependent case:
\[
m=g(t)v_t +v.
\]
We choose $g(t)$ so that no remainder term appears.
From \eqref{eq:p-system}, we obtain
\[
m_t+g(t)(p'(v)v_x)_x=(1+g'(t)-b(t)g(t))v_t.
\]
Thus, we define $g(t)$ by the ordinary differential equation
\begin{equation}\label{eq:g-ode}
1+g'(t)=b(t)g(t).
\end{equation}
With this choice of $g(t)$, we obtain
\begin{equation}\label{eq:m-equation}
m_t+g(t)(p'(v)v_x)_x=0.
\end{equation}
Therefore, the divergence form is derived without differentiability of $b(t)$.

However, the energy method for this problem requires differentiating the equation once with respect to time.
To avoid the appearance of derivatives of $b(t)$ caused by direct differentiation, we introduce another auxiliary function $h(t)$, which allows us to obtain a divergence-form equation after time differentiation.
See Subsection~\ref{sec:reformulations} for details.

Using auxiliary functions $g(t)$ and $h(t)$, we reformulate our problem into two perturbation problems, and derive energy estimates.
As a result, we obtain asymptotic estimates nearly the same as
\eqref{eq:cui-rates} for damping coefficients $b$ in the following class
$\mathcal D$, which does not require differentiability.

\begin{definition}\label{def:damping-class}
(i) Let $b\in C([0,\infty))$ satisfy
\begin{equation}\label{eq:b-bounds}
c_0(1+t)^{-\lambda}\le b(t)\le c_1(1+t)^{-\lambda},
\qquad 0\le\lambda<1,
\qquad c_0,c_1>0.
\end{equation}
We define the minimal positive solution
$g\in C^1([0,\infty))$ of \eqref{eq:g-ode} by
\begin{equation}\label{eq:g-definition}
g(t):=\int_t^\infty
\exp\left(-\int_t^s b(r)\,dr\right)\,ds.
\end{equation}
This is well-defined, which can be verified by using Lemma 3.1 in
\cite{Sobajima'23}.

\noindent(ii) We say that $b\in\mathcal D$ if and only if
$b\in C([0,\infty))$ satisfies \eqref{eq:b-bounds} and the function $g$
defined by \eqref{eq:g-definition} satisfies
\begin{equation}\label{eq:g-prime-bound}
|g'(t)|\le c_2(1+t)^{-(1-\lambda)},
\qquad
c_2>0.
\end{equation}
\end{definition}

The contents of this paper are as follows.  In
Subsection~\ref{sec:reformulations}, we give expected asymptotic profiles, and reformulate our problem
into two perturbation equations around the profiles.
In Subsection~\ref{sec:results}, we state the
results of this paper.
In Subsection~\ref{sec:auxiliary-functions}, we define the two auxiliary functions and state their properties used in this paper.
We also provide a sufficient condition for $b\in\mathcal D$.
In Subsection~\ref{sec:asymptotic-profiles}, we introduce a self-similar diffusion wave, construct the asymptotic profiles using this diffusion wave and the auxiliary functions, and establish estimates for these profiles.  In Section~\ref{sec:energy-estimates}, to prove the main result, we define energy functionals and use them to derive energy estimates for the solutions of the two reformulated problems.

\section{Reformulations and results}\label{sec:reformulations-results}

\subsection{Reformulations}\label{sec:reformulations}
Following~\cite{Nishihara}, we formally reformulate problem \eqref{eq:p-system}--\eqref{eq:end-states} into perturbation equations.

Motivated by \eqref{eq:m-equation} and the expectation that $m_t\approx v_t$ as $t\to\infty$, we take the solution $\bar v$ of the following equation as an asymptotic profile of $v$:
\[
\left\{
\begin{aligned}
&\bar v_t+g(t)(p(\bar v))_{xx}=0,\\
&\bar v(x,t)\to v_\pm \qquad \text{as} \quad x\to\pm\infty.
\end{aligned}
\right.
\]
Moreover, using \eqref{eq:g-ode}, the original system \eqref{eq:p-system} can be rewritten as
\begin{equation}\label{eq:p-system-g}
\left\{
\begin{aligned}
&v_t-u_x=0,\\
&(g(t)u)_t+g(t)(p(v))_x=-u.
\end{aligned}
\right.
\end{equation}
In analogy with Darcy's law, we define the velocity $\bar u$ corresponding to $\bar v$ by
\[
\bar u:=-g(t)(p(\bar v))_x.
\]
Then $(\bar v,\bar u)$ satisfies
\begin{equation}\label{eq:profile-system}
\left\{
\begin{aligned}
&\bar v_t-\bar u_x=0,\\
&g(t)(p(\bar v))_x=-\bar u.
\end{aligned}
\right.
\end{equation}
However, $\bar u(x,t)\to0$ as $x\to\pm\infty$.
We therefore introduce correction functions $(\hat v,\hat u)(x,t)$.
Motivated by the second equation of \eqref{eq:p-system}, we impose
\[
u(x,t)\to u_\pm e^{-B(t)} \qquad \text{as} \quad x\to\pm\infty,
\qquad B(t):=\int_0^t b(r)\,dr.
\]
Let $m_0\in C_c^\infty(\mathbb R)$ satisfy $\int_{-\infty}^{\infty}m_0(y)\,dy=1$, and define
\[
\begin{aligned}
\hat u(x,t)&:=\left[u_-+(u_+-u_-)\int_{-\infty}^x m_0(y)\,dy\right]e^{-B(t)},\\
\hat v(x,t)&:=-(u_+-u_-)m_0(x)g(t)e^{-B(t)}.
\end{aligned}
\]
These functions satisfy
\begin{equation}\label{eq:correction-system}
\left\{
\begin{aligned}
&\hat v_t-\hat u_x=0,\\
&\hat u_t=-b(t)\hat u
\end{aligned}
\right.
\qquad\text{or}\qquad
\left\{
\begin{aligned}
&\hat v_t-\hat u_x=0,\\
&(g(t)\hat u)_t=-\hat u.
\end{aligned}
\right.
\end{equation}
Consequently, after incorporating these correction functions, we have
\[
(u-\bar u-\hat u)(x,t)\to0
\qquad \text{as} \quad x\to\pm\infty.
\]
Furthermore, by choosing a suitable shift $x_0$, we can obtain
\[
\int_{-\infty}^{\infty}
\bigl(v(x,t)-\bar v(x+x_0,t)-\hat v(x,t)\bigr)\,dx=0.
\]
It suffices to take
\[
x_0=\frac{1}{v_+-v_-}\int_{-\infty}^{\infty}
\bigl(v_0(x)-\bar v(x,0)-\hat v(x,0)\bigr)\,dx.
\]
We therefore introduce
\[
V(x,t):=\int_{-\infty}^x
\bigl(v(y,t)-\bar v(y+x_0,t)-\hat v(y,t)\bigr)\,dy
\]
and derive its equation.
Hereafter, we often omit the shift $x_0$ from the notation.
Subtracting \eqref{eq:profile-system} and \eqref{eq:correction-system} from
\eqref{eq:p-system-g}, we obtain
\begin{equation}\label{eq:perturbation-system}
\left\{
\begin{aligned}
&(v-\bar v-\hat v)_t-(u-\bar u-\hat u)_x=0,\\
&\bigl(g(t)(u-\bar u-\hat u)\bigr)_t +g(t)\bigl(p(v)-p(\bar v)\bigr)_x =-(u-\bar u-\hat u)
+\bigl(g(t)^2p(\bar v)_x\bigr)_t.
\end{aligned}
\right.
\end{equation}
Put
\[
\widetilde v(x,t):=\bar v(x+x_0,t)+\hat v(x,t),
\qquad
\widetilde u(x,t):=\bar u(x+x_0,t)+\hat u(x,t).
\]
Since $V_x=v-\widetilde v$ and $V_t=u-\widetilde u$, system \eqref{eq:perturbation-system} leads to the following formulation:
\begin{equation}\label{eq:V-equation}
\left\{
\begin{aligned}
&V_{tt}+\bigl(p'(\widetilde v)V_x\bigr)_x+b(t)V_t=F,
&\qquad& (x,t)\in\mathbb{R}\times(0,\infty),\\
&(V,V_t)(x,0)=(V_0,V_1)(x),
&& x\in\mathbb{R},
\end{aligned}
\right.
\end{equation}
where
\[
F:=\frac{1}{g(t)}\bigl(g(t)^2p(\bar v)_x\bigr)_t
-\bigl(p(V_x+\widetilde v)-p(\bar v)
-p'(\widetilde v)V_x\bigr)_x.
\]
Furthermore, equation \eqref{eq:V-equation} can be rewritten in divergence form:
\begin{equation}\label{eq:V-equation-g}
\bigl(g(t)V_t+V\bigr)_t
+g(t)\bigl(p'(\widetilde v)V_x\bigr)_x=g(t)F.
\end{equation}

To establish a priori estimates, we derive energy estimates for a solution $V$ of \eqref{eq:V-equation}. Under a suitable smallness assumption, the nonlinear term $F$ can be controlled.
More precisely, we derive energy
estimates for $\partial_x^kV(t)$ $(k=0,1,2)$ and
$\partial_x^k\partial_tV(t)$ $(k=0,1)$, where $V$ is the local-in-time solution satisfying
\[
V\in\bigcap_{i=0}^2
C^i\bigl([0,T];H^{3-i}(\mathbb{R})\bigr),
\qquad T>0.
\]
Under the assumption
$b\in C([0,\infty))$ and the smallness assumptions of $(V_0,V_1) \in H^3(\mathbb{R}) \times H^2(\mathbb{R})$, the existence of this local solution follows from
Theorem III of Hughes--Kato--Marsden~\cite{Hughes-Kato-Marsden}.

To obtain the desired energy estimates for $\partial_x^k\partial_t V(t)$, we differentiate \eqref{eq:V-equation} once with respect to time.
However, direct differentiation produces the term $b'(t)$.
Therefore, inspired by \cite{Sobajima'22}, we introduce another auxiliary function $h(t)$ appropriately (as explained below) so that divergence form is restored after time differentiation.
We introduce $W:=V_t/h(t)$ and define $m^*$ by
\[
m^*=h(t)W_t+W.
\]
In terms of $V$, this can be written as
\[
m^*=V_{tt}+\frac{1-h'(t)}{h(t)}V_t.
\]
We choose $h(t)$ to satisfy
\begin{equation}\label{eq:h-ode}
1-h'(t)=b(t)h(t).
\end{equation}
It then follows from \eqref{eq:V-equation} that
\[
m^*+\bigl(p'(\widetilde v)V_x\bigr)_x=F.
\]
Differentiating this identity with respect to time yields divergence form
\begin{equation}\label{eq:Mstar-equation}
m_t^*+h(t)\bigl(p'(\widetilde v)W_x\bigr)_x
=F_t-\bigl(p'(\widetilde v)_tV_x\bigr)_x.
\end{equation}

Although this construction avoids the appearance of $b'(t)$, the term $F_t$ contains $g''(t)$ and thus again requires the differentiability of $b(t)$.
This difficulty arises because $\bar v$ is constructed using $g(t)$.
Therefore, replacing $g(t)$ by $h(t)$ in the construction of $\bar v$, we define a new profile $\bar v^*$ and put
\[
\widetilde v^*(x,t):=\bar v^*(x+x_0^*,t)+\hat v(x,t),
\]
where $x_0^*$ is defined by
\[
x_0^*=\frac{1}{v_+-v_-}\int_{-\infty}^{\infty}
\bigl(v_0(x)-\bar v^*(x,0)-\hat v(x,0)\bigr)\,dx.
\]
Since $\bar v(x,0)=\bar v^*(x,0)$ (see
Subsection~\ref{sec:asymptotic-profiles}), we have
$x_0=x_0^*$.
We then introduce $V^*$ by
\[
V^*(x,t)=\int_{-\infty}^x
\bigl(v-\widetilde v^*\bigr)(y,t)\,dy.
\]
By \eqref{eq:h-ode}, the original system \eqref{eq:p-system} can be rewritten in terms of $h(t)$ as
\[
\left\{
\begin{aligned}
&v_t-u_x=0,\\
&h(t)\left(\frac{u}{h(t)}\right)_t+(p(v))_x
=-\frac{u}{h(t)}.
\end{aligned}
\right.
\]
By an argument analogous to the reformulation for $V$, we obtain
\begin{equation}\label{eq:Vstar-equation}
V_{tt}^*+\bigl(p'(\widetilde v^*)V_x^*\bigr)_x
+b(t)V_t^*=F^*,
\end{equation}
where
\[
F^*:=h(t)(p(\bar v^*))_{xt}
-\bigl(p(V_x^*+\widetilde v^*)-p(\bar v^*)
-p'(\widetilde v^*)V_x^*\bigr)_x.
\]
Then $F^*$ can be differentiated once in time without requiring $b'(t)$.
We redefine $W$ by $W:=V_t^*/h(t)$.
By a calculation analogous to that leading to \eqref{eq:Mstar-equation}, we obtain the following equation:
\begin{equation}\label{eq:Wstar-equation}
(h(t)W_t+W)_t+h(t)\bigl(p'(\widetilde v^*)W_x\bigr)_x
=F_t^*-\bigl(p'(\widetilde v^*)_tV_x^*\bigr)_x.
\end{equation}
To be precise, $V^*$ is defined in terms of $V$ as follows:
\[
V^*(x,t):=V(x,t)+\Theta(x,t), \qquad
\Theta(x,t):=\int_{-\infty}^x
\bigl(\bar v-\bar v^*\bigr)(y+x_0,t)\,dy.
\]
Indeed, this $V^*$ satisfies \eqref{eq:Vstar-equation}.

Equation \eqref{eq:Wstar-equation} has the same structure as \eqref{eq:V-equation-g}.
Therefore, we can use an analogous argument to derive energy estimates.
We use an appropriately modified form of the auxiliary energy introduced
in \cite{Sobajima'25}.

\subsection{Results}\label{sec:results}
Let $b\in\mathcal D$, and put
\[
G(t)=\int_0^t g(r)\,dr,
\qquad
H(t)=\int_0^t h(r)\,dr.
\]
We define
\[
Q_\lambda(t):=
\left\{
\begin{aligned}
&1, &&\qquad 0\le\lambda<1/7,\\
&\left(\log(2+t)\right)^\frac12, &&\qquad \lambda=1/7,\\
&(1+t)^{\frac{-1+7\lambda}{4}}, &&\qquad 1/7<\lambda<1.
\end{aligned}
\right.
\]
Fix a sufficiently large $t_0>0$.
Hereafter, we often write
$\|\cdot\|=\|\cdot\|_{L^2}$.

\begin{proposition}[A priori estimate]\label{prop:apriori}
Assume that $b\in\mathcal D$.
Then there exist constants $\varepsilon>0$ and $C>0$ such that if
\[
\delta:=|v_+-v_-|+|u_+-u_-|\le\varepsilon
\]
holds and $V\in\bigcap_{i=0}^2
C^i\bigl([t_0,t_0+T];H^{3-i}(\mathbb{R})\bigr)$ $(T>0)$
is a solution of \eqref{eq:V-equation} satisfying
\begin{equation}\label{eq:apriori-assumption}
\sup_{t\in[t_0,t_0+T]}Q_\lambda(t)^{-2}
\left(
\sum_{k=0}^3G(t)^k\|\partial_x^kV(t)\|^2
+\sum_{k=0}^2\frac{G(t)^{k+2}}{g(t)^2}
\|\partial_x^k\partial_tV(t)\|^2
\right)
\le \varepsilon^2,
\end{equation}
the following estimate holds:
\begin{equation}\label{eq:decay-estimate}
\begin{aligned}
&\sum_{k=0}^3G(t)^k\|\partial_x^kV(t)\|^2
+\sum_{k=0}^2\frac{G(t)^{k+2}}{g(t)^2}\|\partial_x^k\partial_tV(t)\|^2\\
&\le C\left(\|V(t_0)\|_{H^3}^2+\|V_t(t_0)\|_{H^2}^2+\delta^2 Q_\lambda(t)^2\right),
\qquad t\in[t_0,t_0+T].
\end{aligned}
\end{equation}
\end{proposition}

The proof is given in Subsections~\ref{sec:energy-V} and~\ref{sec:energy-W}.

\begin{remark}
We explain where $Q_\lambda(t)$ in
\eqref{eq:decay-estimate} comes from.
In the energy estimates in Section~\ref{sec:energy-estimates}, for example, we estimate
\[
\int_{t_0}^t g(r)^5(1+G(r))^{-\frac32} \,dr,
\qquad
\int_{t_0}^t h(r)^5(1+H(r))^{-\frac32} \,dr,
\]
which arise from the forcing terms of the equations for $V$ and $W$ in
\eqref{eq:V-equation-g} and \eqref{eq:Wstar-equation}, respectively.
These integrals are bounded by
$Q_\lambda(t)^2$.
In other words, the two terms
$\bigl(g(t)^2p(\bar v)_x\bigr)_t$ and
$\bigl(h(t)p(\bar v^*)_{xt}\bigr)_t$ affect the decay rates.
\end{remark}

By deriving an a priori estimate also on $[0,t_0]$, we construct a
solution on this interval if
$\|V_0\|_{H^3}+\|V_1\|_{H^2}+\delta$ is sufficiently small
(see Subsection~\ref{sec:conclusion}).
Therefore, combining this construction and local existence with Proposition~\ref{prop:apriori}, we
obtain the following theorems.

\begin{theorem}[Global existence]\label{thm:global-existence}
Assume that $b\in\mathcal D$.
For $(V_0,V_1)\in H^3(\mathbb{R})\times H^2(\mathbb{R})$,
if
$\|V_0\|_{H^3}+\|V_1\|_{H^2}+\delta$
is sufficiently small, then \eqref{eq:V-equation} admits a global-in-time solution
$V\in\bigcap_{i=0}^2
C^i\bigl([0,\infty);H^{3-i}\bigr)$.
Moreover, $V(x,t)$ has the estimate \eqref{eq:decay-estimate} for
$t\in [t_0,\infty)$.
\end{theorem}

\begin{theorem}[Main theorem]\label{thm:main}
Assume that $b\in\mathcal D$.
Let $(V_0,V_1)\in H^3(\mathbb{R})\times H^2(\mathbb{R})$, and 
$\|V_0\|_{H^3}+\|V_1\|_{H^2}+\delta$
be sufficiently small.  Let $V(x,t)$ be the global-in-time solution of
\eqref{eq:V-equation} obtained in Theorem~\ref{thm:global-existence}.
Provided $(v_0,u_0)(x)=(V_{0x}(x)+\widetilde v(x,0),
V_1(x)+\widetilde u(x,0))$, the pair
$(v,u):=(V_x+\widetilde v,V_t+\widetilde u)$
is a global solution of \eqref{eq:p-system}--\eqref{eq:end-states} and
satisfies the following estimates:
\[
\begin{aligned}
\|(v-\bar v)(t)\|_{L^2} &= O\bigl((1+G(t))^{-\frac12}Q_\lambda(t)\bigr), \\
\|(u-\bar u-\hat u)(t)\|_{L^2} &= O\bigl(g(t)(1+G(t))^{-1}Q_\lambda(t)\bigr)
\end{aligned}
\]
and
\[
\begin{aligned}
\|(v-\bar v)(t)\|_{L^\infty}&=O\bigl((1+G(t))^{-\frac34}Q_\lambda(t)\bigr), \\
\|(u-\bar u)(t)\|_{L^\infty}&=O\bigl(g(t)(1+G(t))^{-\frac54}Q_\lambda(t)\bigr)
\end{aligned}
\]
as $t\to\infty$.
Namely,
\[
\|(v-\bar v)(t)\|_{L^\infty}=
\left\{
\begin{aligned}
&O\bigl((1+t)^{-\frac{3(1+\lambda)}{4}}\bigr),
&\qquad&0\le\lambda<1/7,\\
&O\bigl((1+t)^{-\frac67}(\log(2+t))^{\frac12}\bigr),
&&\lambda=1/7,\\
&O\bigl((1+t)^{-(1-\lambda)}\bigr),
&&1/7<\lambda<1,
\end{aligned}
\right.
\]
and
\[
\|(u-\bar u)(t)\|_{L^\infty}=
\left\{
\begin{aligned}
&O\bigl((1+t)^{-\frac{5+\lambda}{4}}\bigr),
&\qquad&0\le\lambda<1/7,\\
&O\bigl((1+t)^{-\frac97}(\log(2+t))^{\frac12}\bigr),
&&\lambda=1/7,\\
&O\bigl((1+t)^{-\frac{3(1-\lambda)}{2}}\bigr),
&&1/7<\lambda<1
\end{aligned}
\right.
\]
as $t\to\infty$.
Here $(\bar v,\bar u)$ denotes the shifted profile
$(\bar v,\bar u)(x+x_0,t)$.
\end{theorem}

\begin{remark}
The $L^\infty$ estimates mean both $v-\bar v$ and
$u-\bar u$ decay for every $0\le\lambda<1$.
On the other hand, when $\lambda$ is close to $1$, the $L^2$ estimates do not mean they decay.  More precisely, $v-\bar v$ decays for $0\le\lambda<3/5$, while $u-\bar u-\hat u$ decays for
$0\le\lambda<5/7$.
\end{remark}

\begin{remark}
Suppose that $b(t)=(1+t)^{-\lambda}$.  The $L^\infty$ asymptotic estimates of
Theorem~\ref{thm:main} remain valid when $(\bar v,\bar u)$ is replaced
by the profile introduced in \cite{Cui} (see
Remark~\ref{rem:profile-comparison} below).
\end{remark}

\begin{remark}
For $\lambda\ne 1/7$, the $L^\infty$ decay rates in Theorem~\ref{thm:main} coincide with those obtained by Cui~\cite{Cui} in \eqref{eq:cui-rates}.
In the critical case $\lambda=1/7$, the polynomial loss in \eqref{eq:cui-rates} is improved to a logarithmic loss.
\end{remark}

\begin{remark}\label{rem:nonsmooth-damping}
The class $\mathcal D$ contains certain non-differentiable or rapidly
oscillating damping coefficients.
For example, for $|\nu|<1$ and $k>0$,
\[
b(t)=(1+t)^{-\lambda}\bigl(1+\nu\sin(\log(2+t))\bigr)
\left(1+\frac{|\sin(t^k)|}{(1+t)^{1-\lambda}}\right)
\]
belongs to $\mathcal D$.  Moreover, this example shows that the condition
$b\in\mathcal D$ does not require the existence of
$\lim_{t\to\infty}b(t)(1+t)^\lambda$.  See
Lemma~\ref{lem:admissible-damping} below.
\end{remark}

\section{Auxiliary functions and asymptotic profiles}\label{sec:preliminaries}

\subsection{Auxiliary functions}\label{sec:auxiliary-functions}

We first introduce two auxiliary functions $g(t)$ and $h(t)$.
If $b\in\mathcal D$,
then $g(t)$ satisfies
\[
\lim_{t\to\infty}g'(t)=0,
\qquad
\lim_{t\to\infty}b(t)g(t)=1
\]
(see Definition~\ref{def:damping-class}).
Thus, the damping coefficients $b(t)$ considered in this paper belong to the class
treated in \cite{Sobajima'23, Sobajima'25}.

We define $h\in C^1([0,\infty))$ as a solution of \eqref{eq:h-ode}.  Solving
\eqref{eq:h-ode}, we obtain
\begin{equation}\label{eq:h-formula}
h(t)=\exp\left(-\int_0^t b(r)\,dr\right)
\left[
h(0)+\int_0^t\exp\left(\int_0^s b(r)\,dr\right)\,ds
\right].
\end{equation}
We fix an arbitrary $h(0)>0$ and define the positive function
$h(t)$ by \eqref{eq:h-formula}.
Then $h'(t)$ has the same decay upper bound as $g'(t)$
(see Lemma~\ref{lem:h-prime} below).

\medskip
\noindent\textbf{Notation.}
For nonnegative functions $p(t)$ and $q(t)$, we write
$p(t)\lesssim q(t)$ if $p(t)\le cq(t)$ holds for some constant $c>0$
independent of $t$.  We write $p(t)\asymp q(t)$ if
$p(t)\lesssim q(t)$ and $q(t)\lesssim p(t)$ hold.

\begin{lemma}\label{lem:h-prime}
Suppose that $b\in\mathcal D$.  Then there exists a constant $C_3>0$ such that
\[
|h'(t)|\le C_3(1+t)^{-(1-\lambda)},
\qquad t\ge0.
\]
\end{lemma}

\begin{proof}
Adding \eqref{eq:g-ode} and \eqref{eq:h-ode}, we obtain
\[
2+g'(t)-h'(t)=b(t)g(t)+b(t)h(t).
\]
Put $\theta(t):=g(t)-h(t)$.  Using \eqref{eq:g-ode} again, we obtain
\[
\theta'(t)+b(t)\theta(t)=2g'(t).
\]
Solving this equation, we obtain
\[
\theta(t)=\theta(0)\exp\left(-\int_0^t b(r)\,dr\right)
+2\int_0^t g'(s)
\exp\left(-\int_s^t b(r)\,dr\right)\,ds.
\]
To estimate the last term,
we divide the integral interval into
$[0,t/2]$ and $[t/2,t]$:
\[
\begin{aligned}
\left|\int_0^{t/2}g'(s)
\exp\left(-\int_s^t b(r)\,dr\right)\,ds\right|
&\lesssim \exp\left(-\int_{t/2}^t b(r)\,dr\right)
\int_0^{t/2}(1+s)^{-(1-\lambda)}\,ds, \\
\left|\int_{t/2}^t g'(s)
\exp\left(-\int_s^t b(r)\,dr\right)\,ds\right|
&=\left|\int_{t/2}^t \frac{g'(s)}{b(s)}
\frac{d}{ds}\left[\exp\left(-\int_s^t b(r)\,dr\right)\right]\,ds\right|\\
&\le \sup_{s\in[t/2,t]}\left|\frac{g'(s)}{b(s)}\right|
\left|1-\exp\left(-\int_{t/2}^t b(r)\,dr\right)\right|\\
&\lesssim(1+t)^{-(1-2\lambda)}.
\end{aligned}
\]
The integral over $[0,t/2]$ decays exponentially and hence does not affect
the polynomial estimate for the integral over $[t/2,t]$.  Consequently,
\[
|\theta(t)|\lesssim (1+t)^{-(1-2\lambda)}.
\]
Since
$h'(t)=-g'(t)+b(t)\theta(t)$,
we obtain
\[
|h'(t)|\lesssim (1+t)^{-(1-\lambda)}.
\]
\end{proof}

By Lemma~\ref{lem:h-prime}, $h(t)$ also satisfies
\[
\lim_{t\to\infty}h'(t)=0,
\qquad
\lim_{t\to\infty}b(t)h(t)=1.
\]
Hence, for $b\in\mathcal D$, the auxiliary functions satisfy
\begin{subequations}\label{eq:auxiliary-orders}
\begin{equation}\label{eq:gh-orders}
g(t)\asymp(1+t)^\lambda,
\qquad
h(t)\asymp(1+t)^\lambda,
\end{equation}
and
\begin{equation}\label{eq:gh-prime-orders}
|g'(t)|+|h'(t)|\lesssim(1+t)^{-(1-\lambda)}.
\end{equation}
Moreover, \eqref{eq:gh-orders} implies
\begin{equation}\label{eq:GH-orders}
1+G(t)\asymp(1+t)^{1+\lambda},
\qquad
1+H(t)\asymp(1+t)^{1+\lambda},
\end{equation}
\end{subequations}
where
\[
G(t):=\int_0^t g(r)\,dr,
\qquad
H(t):=\int_0^t h(r)\,dr.
\]

We now give a sufficient condition for $b \in\mathcal D$.
The following lemma shows that the example in
Remark~\ref{rem:nonsmooth-damping} belongs to $\mathcal D$.

\begin{lemma}\label{lem:admissible-damping}
Let $b_0\in C^1([0,\infty))$ satisfy
\[
b_0(t)\asymp(1+t)^{-\lambda},
\qquad
|b_0'(t)|\lesssim(1+t)^{-\lambda-1},
\]
and let $\widetilde b\in C([0,\infty))$ satisfy
$1+\widetilde b(t)>0$ and
\[
|\widetilde b(t)|\lesssim(1+t)^{-(1-\lambda)}.
\]
Then
\[
b(t):=b_0(t)\bigl(1+\widetilde b(t)\bigr)
\]
belongs to $\mathcal D$.
\end{lemma}

\begin{proof}
The function $b(t)=b_0(t)(1+\widetilde b(t))$ satisfies
\eqref{eq:b-bounds}.  We therefore define $g(t)$ by
\eqref{eq:g-definition}.

We first show that $b(t)g(t)\lesssim1$.  Put
\[
b_*(t):=c_0(1+t)^{-\lambda}\le b(t)
\]
and define
\[
g_*(t):=\int_t^\infty
\exp\left(-\int_t^s b_*(r)\,dr\right)\,ds.
\]
Then \eqref{eq:g-definition} gives $g(t)\le g_*(t)$.  Moreover, by
Lemma 2.1 of \cite{Sobajima'22},
\[
\lim_{t\to\infty}b_*(t)g_*(t)=1.
\]
Since \eqref{eq:b-bounds} gives
$b(t)\le (c_1/c_0) b_*(t)$, we obtain
\begin{equation}\label{eq:bg-bound}
b(t)g(t)\le\frac{c_1}{c_0}b_*(t)g_*(t)\lesssim1.
\end{equation}

For $b_0$, we also define an auxiliary function
$g_0\in C^1([0,\infty))$ as the minimal positive solution of
\begin{equation}\label{eq:g0-ode}
1+g_0'(t)=b_0(t)g_0(t).
\end{equation}
By the same argument as in Lemma 2.2 of \cite{Sobajima'22},
$g_0$ satisfies
\begin{equation}\label{eq:g0-prime-bound}
g_0(t)\asymp(1+t)^\lambda,
\qquad
|g_0'(t)|\lesssim(1+t)^{-(1-\lambda)}.
\end{equation}

Subtracting \eqref{eq:g0-ode} from \eqref{eq:g-ode}, we obtain
\[
g'(t)-g_0'(t)
=b(t)\bigl(g(t)-g_0(t)\bigr)
+b_0(t)g_0(t)\widetilde b(t).
\]
Put $\kappa(t):=g(t)-g_0(t)$.  Then
\[
\kappa'(t)-b(t)\kappa(t)=b_0(t)g_0(t)\widetilde b(t).
\]
Solving this equation, we obtain
\[
\kappa(t)=-\int_t^\infty b_0(s)g_0(s)\widetilde b(s)
\exp\left(-\int_t^s b(r)\,dr\right)\,ds.
\]
Note that the boundary term at infinity vanishes.
Using $b_0(t)g_0(t)\lesssim1$ and the assumption on $\widetilde b$, we have
\[
|\kappa(t)|\lesssim(1+t)^{-(1-\lambda)}
\int_t^\infty \exp\left(-\int_t^s b(r)\,dr\right)\,ds
=(1+t)^{-(1-\lambda)}g(t).
\]
Combining the identity
\[
g'(t)=g_0'(t)+b(t)\kappa(t)
+b_0(t)g_0(t)\widetilde b(t)
\]
with \eqref{eq:bg-bound} and
\eqref{eq:g0-prime-bound}, we obtain
\eqref{eq:g-prime-bound}.
\end{proof}

\subsection{Asymptotic profiles}\label{sec:asymptotic-profiles}

For the problem \eqref{eq:p-system}--\eqref{eq:end-states} with constant
damping $b(t)\equiv\alpha>0$, a self-similar solution satisfying the
following nonlinear diffusion equation was used in
\cite{Hsiao-Liu,Nishihara} as an asymptotic profile of $v(x,t)$: 
\begin{equation}\label{eq:standard-diffusion}
\left\{
\begin{aligned}
&\bar v_t +\frac{1}{\alpha}(p(\bar v))_{xx}=0,
    \qquad (x,t)\in \mathbb{R}\times(0,\infty), \\
&\bar v(x,t)\to v_{\pm} \qquad \text{as} \quad x\to\pm\infty.
\end{aligned}
\right.
\end{equation}

We instead
introduce two asymptotic profiles satisfying
\begin{equation}\label{eq:time-dependent-diffusion}
\begin{aligned}
\bar v_t+g(t)(p(\bar v))_{xx}&=0, \\
\bar v_t^*+h(t)(p(\bar v^*))_{xx}&=0,
\end{aligned}
\end{equation}
and the same far-field condition as in \eqref{eq:standard-diffusion}.
They are constructed by rescaling the time variable in the
self-similar solution of \eqref{eq:standard-diffusion}.
We now describe this construction.

First, we consider a self-similar
solution of
\begin{equation}\label{eq:sigma}
\sigma_t+(p(\sigma))_{xx}=0.
\end{equation}
Substituting
\[
\sigma(x,t)=\phi(\xi), \qquad
\xi:=\frac{x}{\sqrt{1+t}},
\]
into \eqref{eq:sigma} yields the ordinary differential equation
\begin{equation}\label{eq:phi-ode}
-\bigl(p'(\phi)\phi'\bigr)'+\frac12\xi\phi'=0,
\qquad
\xi\in\mathbb{R}.
\end{equation}
By van Duyn--Peletier~\cite{Duyn}, together with the
assumption $p\in C^4((0,\infty))$, equation~\eqref{eq:phi-ode} admits a unique solution
$\phi\in C^4(\mathbb{R})$ satisfying
\[
\lim_{\xi\to\pm\infty}\phi(\xi)=v_{\pm}.
\]
Moreover, $\phi$ is strictly monotone.  An argument analogous to that in
Section~2 of \cite{Hsiao-Liu} shows that the derivatives of $\phi$ satisfy the
Gaussian decay estimate
\begin{equation}\label{eq:phi-gaussian}
\sum_{k=1}^4|\phi^{(k)}(\xi)|
\le C|v_+-v_-|e^{-c\xi^2}
\end{equation}
for some positive constants $C$ and $c$.

Using $\phi$, we define
\begin{equation}\label{eq:profile-definitions}
\bar v(x,t):=\phi\left(\frac{x}{\sqrt{1+G(t)}}\right), \qquad
\bar v^*(x,t):=\phi\left(\frac{x}{\sqrt{1+H(t)}}\right).
\end{equation}
These functions satisfy \eqref{eq:time-dependent-diffusion} and the
far-field conditions
\[
\bar v(x,t)\to v_{\pm}, \qquad
\bar v^*(x,t)\to v_{\pm}
\qquad \text{as} \quad x\to\pm\infty.
\]
\eqref{eq:phi-gaussian} gives the following estimates.

\begin{lemma}\label{lem:profile-decay}
Assume that $p \in C^4((0,\infty))$.
Then the profile $\bar v$ (defined in \eqref{eq:profile-definitions}) satisfies
\[
\begin{aligned}
\|\partial_x^k \bar v\|_{L^2}
&\lesssim \delta(1+G(t))^{\frac14-\frac k2},
&\qquad& k=1,2,3,4,\\
\|\partial_x^k\partial_t \bar v\|_{L^2}
&\lesssim \delta g(t)(1+G(t))^{-\frac34-\frac k2},
&& k=0,1,2,3,\\
\|\partial_x^k\partial_t^2 \bar v\|_{L^2}
&\lesssim \delta g(t)^2(1+G(t))^{-\frac74-\frac k2},
&& k=0,1,2.
\end{aligned}
\]
The corresponding estimates for $\bar v^*$ are obtained by replacing $g(t)$
and $G(t)$ with $h(t)$ and $H(t)$, respectively.
\end{lemma}

We note that the third estimate uses the following bounds:
\begin{equation}\label{eq:g2-over-G}
|g'(t)|\lesssim \frac{g(t)^2}{1+G(t)},
\qquad
|h'(t)|\lesssim \frac{h(t)^2}{1+H(t)}.
\end{equation}

We next define
\[
\bar u(x,t):=-g(t)(p(\bar v(x,t)))_x, \qquad
\bar u^*(x,t):=-h(t)(p(\bar v^*(x,t)))_x.
\]
Then $(\bar v,\bar u)$ satisfies \eqref{eq:profile-system} and
\begin{equation}\label{eq:profile-far-field}
(\bar v,\bar u)(x,t) \to (v_\pm,0) \qquad
\text{as} \quad x\to\pm\infty.
\end{equation}
Likewise, $(\bar v^*,\bar u^*)$ satisfies the equation corresponding to \eqref{eq:profile-system}
 (replace $g(t)$ with $h(t)$) and the same far-field condition as in \eqref{eq:profile-far-field}.

\begin{remark}\label{rem:profile-comparison}
When $b(t)=(1+t)^{-\lambda}$, the profile $(\bar v_b,\bar u_b)$ introduced
in \cite{Cui} is given by
\[
\bar v_b(x,t)
:=\phi\left(\frac{x}{\sqrt{(1+t)^{1+\lambda}/(1+\lambda)}}\right), \qquad
\bar u_b(x,t):=-\frac{1}{b(t)}(p(\bar v_b(x,t)))_x,
\]
which satisfies \eqref{eq:variable-profile}.
By \eqref{eq:phi-gaussian}, we obtain
\[
\|\bar v-\bar v_b\|_{L^2}
\lesssim \delta(1+t)^{\frac{-3+5\lambda}{4}},
\qquad
\|\bar v-\bar v_b\|_{L^\infty}
\lesssim \delta(1+t)^{-(1-\lambda)},
\]
and
\[
\|\bar u-\bar u_b\|_{L^2}
\lesssim \delta(1+t)^{\frac{-5+7\lambda}{4}},
\qquad
\|\bar u-\bar u_b\|_{L^\infty}
\lesssim \delta(1+t)^{-\frac{3(1-\lambda)}{2}}.
\]
Thus, Theorem~\ref{thm:main} remains valid when $(\bar v,\bar u)$ is replaced
by $(\bar v_b,\bar u_b)$.
\end{remark}

Here, as a preliminary, we state estimates for the difference between
the two profiles.
\begin{lemma}\label{lem:tilde-difference}
For $\Theta(x,t)=\int_{-\infty}^x (\bar v -\bar v^*)(y+x_0,t) \,dy$, the following estimates hold for any $t\ge 0$:
\[
\begin{aligned}
\| \partial_x^k \Theta(t) \|
&\lesssim \delta (1+t)^{\frac{-1+7\lambda}{4}-\frac{1+\lambda}{2}k},
&\qquad& k=0,1,2,3, \\
\|\partial_x^k\partial_t\Theta(t)\|
&\lesssim \delta
(1+t)^{\frac{-5+7\lambda}{4}-\frac{1+\lambda}{2}k},
&& k=0,1,2.
\end{aligned}
\]
\end{lemma}

\begin{proof}
We first estimate $\|\Theta(t)\|$.
Put $\sigma(x,r)=\phi(x/\sqrt{1+r})$.
Then $\sigma$ satisfies $\sigma_r=-p(\sigma)_{xx}$.
Using the definitions of $\bar v$ and $\bar v^*$ and \ref{eq:phi-gaussian}, we have
\[
\begin{aligned}
\Theta(x,t) &= \int_{-\infty}^x \int_{H(t)}^{G(t)} \sigma_r(y+x_0,r)\,drdy \\
&=-\int_{-\infty}^x \int_{H(t)}^{G(t)} p(\sigma(y+x_0,r))_{yy}\,drdy \\
&=-\int_{H(t)}^{G(t)} p(\sigma(x+x_0,r))_x dr.
\end{aligned}
\]
Moreover, since
$\|p(\sigma(r))_x\| \lesssim \delta(1+r)^{-\frac14}$,
it follows that
\[
\begin{aligned}
\|\Theta(t)\| &\le \int_{\min\{G(t),H(t)\}}^{\max\{G(t),H(t)\}} \|p(\sigma(r))_x\| \, dr \\
&\lesssim \delta (1+\min\{G(t),H(t)\})^{-\frac14}|G(t)-H(t)| \\
&\lesssim \delta (1+t)^{-\frac{1+\lambda}{4}}(1+t)^{2\lambda}.
\end{aligned}
\]
Here we used $G(t)-H(t)=g(t)h(t)-g(0)h(0)$, which follows from $(g(t)h(t))'=g(t)-h(t)$ based on \eqref{eq:g-ode} and \eqref{eq:h-ode}.
This completes the proof of the first estimate for $k=0$.
The cases $k=1,2,3$ are proved similarly.

Next, we estimate $\|\Theta_t(t)\|$.
We have
\[
\begin{aligned}
\Theta_t(t)&=-\Big[g(t)p\big(\sigma(G(t))\big)_x-h(t)p\big(\sigma(H(t))\big)_x\Big] \\
&=-\big[g(t)-h(t)\big]p\big(\sigma(G(t))\big)_x -h(t)\Big[p\big(\sigma(G(t))\big)_x-p\big(\sigma(H(t))\big)_x \Big].
\end{aligned}
\]
By using $\|p\big(\sigma(G(t))\big)_x\|\lesssim \delta(1+G(t))^{-\frac14}$, $\|p\big(\sigma(t)\big)_{xt}\|\lesssim \delta(1+t)^{-\frac54}$ and mean value theorem, it follows that
\[
\begin{aligned}
\|\Theta_t(t)\|
&\lesssim \delta \left|(g(t)h(t))'\right|(1+G(t))^{-\frac14} \\
&\quad+ \delta h(t)\left|G(t)-H(t)\right| (1+\min\{G(t),H(t)\})^{-\frac54} \\
&\lesssim \delta(1+t)^{-1+2\lambda}(1+t)^{-\frac{1+\lambda}{4}}
+\delta(1+t)^\lambda (1+t)^{2\lambda} (1+t)^{-\frac{5(1+\lambda)}{4}}.
\end{aligned}
\]
The cases $k=1,2$ are proved similarly.
\end{proof}

\section{Energy estimates}\label{sec:energy-estimates}

In Subsections~\ref{sec:energy-V} and~\ref{sec:energy-W}, we derive energy
estimates under the assumptions of Proposition~\ref{prop:apriori} in order
to prove it.
In Subsection~\ref{sec:conclusion}, we derive energy estimates on $[0,t_0]$
to ensure the existence of a solution on this interval.
Hereafter, $C$ denotes a generic positive constant
which may vary from line to line.
The domain of integration $\mathbb{R}$ is omitted unless confusion arises.

\subsection{Energy estimates for $V$}\label{sec:energy-V}

In this subsection, we successively define energy functionals satisfying
\[
\begin{aligned}
E_k(t)&\asymp \|\partial_x^{k+1} V(t)\|^2+\|\partial_x^k V_t(t)\|^2, \\
\widetilde E_k(t)&\asymp g(t)^2E_k(t)+\|\partial_x^kV(t)\|^2, \qquad k=0,1,2
\end{aligned}
\]
and derive the corresponding estimates. All estimates obtained in this subsection are summarized in a table at the end.

We note that the a priori assumption
\eqref{eq:apriori-assumption} and Sobolev's inequality yield
\begin{equation}\label{eq:V-linfty}
\begin{aligned}
\|\partial_x^kV(t)\|_{L^\infty}
&\le C\varepsilon 
(1+G(t))^{-\frac14-\frac{k}{2}}Q_\lambda(t),
&\qquad& k=0,1,2, \\
\|\partial_x^k\partial_tV(t)\|_{L^\infty}
&\le C\varepsilon g(t)
(1+G(t))^{-\frac54-\frac{k}{2}}Q_\lambda(t),
&& k=0,1.
\end{aligned}
\end{equation}
Since $V^*=V+\Theta$, Lemma~\ref{lem:tilde-difference} shows that
\eqref{eq:apriori-assumption} (and hence \eqref{eq:V-linfty}) remains valid
with $V$ replaced by $V^*$, up to multiplicative constants.
We shall frequently use \eqref{eq:auxiliary-orders} and
\eqref{eq:g2-over-G}.

We decompose $F^*$ as $F^*=F_1^*+F_2^*(V^*)+F_3^*$,
where
\[
\begin{aligned}
F_1^*&:=h(t)(p(\bar v^*))_{xt},\\
F_2^*(V^*)&:=-\bigl(p(V_x^*+\widetilde v^*)-p(\widetilde v^*)
-p'(\widetilde v^*)V_x^*\bigr)_x,\\
F_3^*&:=-\bigl(p(\widetilde v^*)-p(\bar v^*)\bigr)_x.
\end{aligned}
\]
Using \eqref{eq:g-ode}, equation~\eqref{eq:Vstar-equation} can be rewritten as
\begin{equation}\label{eq:Vstar-g}
\bigl(g(t)V_t^*+V^*\bigr)_t
+g(t)\bigl(p'(\widetilde v^*)V_x^*\bigr)_x=g(t)F^*.
\end{equation}

We define the following energies, which are suitably modified versions of those introduced in~\cite{Sobajima'25}.

\begin{definition}\label{def:E0}
We define the auxiliary energies and usual energies by
\[
\begin{aligned}
\widetilde E_0(t)
&:=-g(t)^2\int p'(\widetilde v)V_x^2\,dx
+\left\|g(t)V_t+\frac12V\right\|^2
+\frac14\|V\|^2,\\
E_0(t)
&:=-\int p'(\widetilde v)V_x^2\,dx+\|V_t\|^2,
\end{aligned}
\]
respectively.
We also define $\widetilde E_0^*(t)$ and $E_0^*(t)$ by replacing
$\widetilde v$ and $V$ with $\widetilde v^*$ and $V^*$, respectively,
in the above definitions.
\end{definition}

\begin{lemma}\label{lem:E0}
If $\varepsilon>0$ is sufficiently small, then the following estimate holds
for $t\in [t_0,t_0+T]$:
\[
\widetilde E_0^*(t)+\int_{t_0}^t g(r)E_0^*(r)\,dr
\le C\left(\|V^*(t_0)\|_{H^1}^2+\|V^*_t(t_0)\|_{L^2}^2
+\delta^2Q_\lambda(t)^2\right).
\]
\end{lemma}

\begin{proof}

Differentiating $\widetilde E_0^*(t)$ and using \eqref{eq:Vstar-g}, we obtain
\[
\begin{aligned}
\frac{d}{dt}\widetilde E_0^*
&=-2gg'\int p'(\widetilde v^*){V_x^*}^2\,dx-g^2\int p''(\widetilde v^*)\widetilde v_t^*{V_x^*}^2\,dx
-2g^2\int p'(\widetilde v^*)V_x^*V_{xt}^*\,dx\\
&\quad+2\int\left(gV_t^*+\frac12V^*\right)
\left(-g\bigl(p'(\widetilde v^*)V_x^*\bigr)_x
-\frac12V_t^*+gF^*\right)\,dx\\
&\quad+\frac12\int V^*V_t^*\,dx\\
&=-2gg'\int p'(\widetilde v^*){V_x^*}^2\,dx
-g^2\int p''(\widetilde v^*)\widetilde v_t^*{V_x^*}^2\,dx\\
&\quad-gE_0^*
+2g\int F^*\left(gV_t^*+\frac12V^*\right)\,dx.
\end{aligned}
\]
Here, we used integration by parts.
We take $t_0\ge0$ sufficiently large so that
$|g'(t)|\le1/4$ for $t\ge t_0$.
Moreover, we have
\[
\|\widetilde v_t^*\|_{L^\infty}
\le C\delta h(1+H)^{-1},
\]
and the same estimate holds with $h$ and $H$ replaced by $g$ and $G$,
respectively.
Then it follows that
\[
\frac{d}{dt}\widetilde E_0^*
\le-\left(\frac12-C\delta\right)gE_0^*
+2g\int F^*\left(gV_t^*+\frac12V^*\right)\,dx.
\]
Here, we used \eqref{eq:g2-over-G} and $\|V^*_x(t)\|^2 \le C E^*_0(t)$.

We next estimate the terms containing $F^*$.  By Young's inequality and
Lemma~\ref{lem:profile-decay}, for $\mu> 0$, we have
\[
\begin{aligned}
g^2\left|\int F_1^*V_t^*\,dx\right|
&\le\mu g\|V_t^*\|^2
+C_\mu\delta^2g^3h^4(1+H)^{-\frac52},\\
g\left|\int F_1^*V^*\,dx\right|
&=gh\left|-\int p(\bar v^*)_tV_x^*\,dx\right|\le\mu g\|V_x^*\|^2
+C_\mu\delta^2gh^4(1+H)^{-\frac32}.
\end{aligned}
\]
For $F_3^*$, we have
\[
\begin{aligned}
2g\left|\int
\bigl(p(\widetilde v^*)-p(\bar v^*)\bigr)_x
\left(gV^*_t+\frac12V^*\right)\,dx\right|
&\le\mu g\bigl(\|V^*_t\|^2+\|V^*_x\|^2\bigr)
+C_\mu\delta^2e^{-B}.
\end{aligned}
\]
Among these terms, $gh^4(1+H)^{-\frac32}$ has the slowest decay.
Consequently,
\[
\frac{d}{dt}\widetilde E_0^*
\le-\left(\frac12-C(\delta+\mu)\right)gE_0^*
+C_\mu\delta^2gh^4(1+H)^{-\frac32}
+2g\|F^*_2(V^*)\|(\widetilde E_0^*)^{\frac12},
\]
on $[t_0,t_0+T]$.
Take $\delta$ and $\mu$ small so that $1/2 - C(\delta + \mu)>0$.
Then the Gronwall-type inequality (Lemma~2.1 of \cite{Sobajima'25}) yields
\[
\begin{aligned}
&\widetilde E_0^*(t) + \left(\frac12-C(\delta+\mu)\right)\int_{t_0}^t g(r)E_0^*(r)\,dr \\
&\le
\left[ \left(\widetilde E_0^*(t_0) + C_\mu\delta^2\int_{t_0}^t g(r)h(r)^4(1+H(r))^{-\frac32} \, dr \right)^\frac12 +\int_{t_0}^t g(r)\|F^*_2(V^*)(r)\| \, dr \right]^2.
\end{aligned}
\]
Here, since
$g(t)h(t)^4(1+H(t))^{-\frac32}
\lesssim(1+t)^{\frac{-3+7\lambda}{2}}$,
its integral produces $Q_\lambda(t)^2$.
Moreover, by the definition of $\widetilde E_0^*(t)$, we have
\[
\widetilde E_0^*(t_0)
\le C\left(\|V^*(t_0)\|_{H^1}^2
+\|V_t^*(t_0)\|_{L^2}^2\right).
\]
On the other hand, $F^*_2(V^*)$ can be written as
\[
F^*_2(V^*)
=-\bigl(p'(V^*_x+\widetilde v^*)-p'(\widetilde v^*)\bigr)V^*_{xx}-\bigl(p'(V^*_x+\widetilde v^*)-p'(\widetilde v^*)
-p''(\widetilde v^*)V^*_x\bigr)\widetilde v^*_x,
\]
and hence
\[
|F^*_2(V^*)|
\le C\bigl(|V^*_{xx}|+|\widetilde v^*_x||V^*_x|\bigr)|V^*_x|.
\]
Therefore, using
$\|\widetilde v^*_x\|_{L^\infty}
\le C\delta(1+H)^{-\frac12}$
and $V^*$ version of \eqref{eq:V-linfty}, we obtain
\[
g\|F^*_2(V^*)\|
\le C(\varepsilon+\delta\varepsilon)
g(1+G)^{-\frac54}Q_\lambda\|V^*_x\|.
\]
By the Cauchy--Schwarz inequality,
\[
\begin{aligned}
&\left(\int_{t_0}^t g(r)\|F^*_2(V^*)(r)\|\,dr\right)^2 \\
&\le C(\varepsilon+\delta\varepsilon)^2
\left(\int_{t_0}^t g(r)(1+G(r))^{-\frac52}Q_\lambda(r)^2\,dr\right)
\left(\int_{t_0}^t g(r)\|V^*_x(r)\|^2\,dr\right)\\
&\le C(\delta+\varepsilon)
\int_{t_0}^t g(r)||V^*_x(r)||^2\,dr.
\end{aligned}
\]
Here,
$g(t)(1+G(t))^{-\frac52}Q_\lambda(t)^2$
decays faster than $(1+t)^{-1}$ for every $0\le\lambda<1$.
Since $0\le\delta(\le\varepsilon)$ and $\varepsilon>0$ is sufficiently
small, we will bound higher powers of $\delta$ and $\varepsilon$
by $\delta$ and $\varepsilon$, respectively.
Therefore, taking $\varepsilon>0$ sufficiently small, we obtain the desired estimate.
\end{proof}

Since
\[
\widetilde E_0(t)\asymp g(t)^2E_0(t)+\|V(t)\|^2,
\qquad
E_0(t)\asymp\|V_x(t)\|^2+\|V_t(t)\|^2,
\]
Lemma~\ref{lem:tilde-difference} and Lemma~\ref{lem:E0} yield the following
estimate.

\begin{corollary}\label{cor:E0}
If $\varepsilon>0$ is sufficiently small, then the following estimate holds
for $t \in [t_0, t_0+T]$:
\[
\widetilde E_0(t)+\int_{t_0}^t g(r)E_0(r)\,dr
\le C\left(\|V(t_0)\|_{H^1}^2+\|V_t(t_0)\|_{L^2}^2
+\delta^2Q_\lambda(t)^2\right).
\]
\end{corollary}

\begin{remark}\label{rem:Vstar-energy}
Here, we explain why we use $V^*$ in Lemma~\ref{lem:E0} instead of $V$.
If the same calculation is carried out for $V$, an estimate of the following
term is required:
\[
g\left|\int F_1V\,dx\right|,
\qquad
F_1:=\frac{1}{g}\bigl(g^2p(\bar v)_x\bigr)_t.
\]
Hence, integration by parts in this term produces
\[
2g|g'|\left|\int p(\bar v)V_x\,dx\right|,
\]
which worsens the decay rate.  On the other hand, if the term is left as $2g\|F_1\|(\widetilde E_0)^{\frac12}$ without integration by parts, the
logarithmic loss for $\lambda=1/7$ becomes $(\log(2+t))^2$.
In the present proof, since $F_1^*=h(p(\bar v^*))_{xt}$, this term does not
appear.
\end{remark}

We also decompose $F$ as $F=F_1+F_2(V)+F_3$,
where
\[
\begin{aligned}
F_1&:=\frac{1}{g(t)}(g(t)^2p(\bar v)_x)_t,\\
F_2(V)&:=-\bigl(p(V_x+\widetilde v)-p(\widetilde v)
-p'(\widetilde v)V_x\bigr)_x,\\
F_3&:=-\bigl(p(\widetilde v)-p(\bar v)\bigr)_x.
\end{aligned}
\]

\begin{lemma}\label{lem:E0-weighted}
If $\varepsilon>0$ is sufficiently small, then the following estimate holds
on $[t_0,t_0+T]$:
\[
GE_0+\int_{t_0}^t bG\|V_t\|^2\,dr
\le C\left(\|V(t_0)\|_{H^1}^2+\|V_t(t_0)\|_{L^2}^2
+\delta^2Q_\lambda^2\right).
\]
\end{lemma}

\begin{proof}
Differentiating $G(t)E_0(t)$ and using \eqref{eq:V-equation}, we obtain
\[
\begin{aligned}
\frac{d}{dt}\bigl(GE_0\bigr)
=&gE_0+G\Biggl(-\int p''(\widetilde v)
\widetilde v_tV_x^2\,dx
-2\int p'(\widetilde v)V_xV_{xt}\,dx\\
&\qquad+2\int V_t\left(-\bigl(p'(\widetilde v)V_x\bigr)_x
-bV_t+F\right)dx\Biggr)\\
\le& CgE_0-2bG\|V_t\|^2
+2G\int FV_t\,dx.
\end{aligned}
\]
The terms involving $F_1$ and $F_3$ satisfy
\[
G\left|\int (F_1+F_3)V_t\,dx\right|
\le \mu bG\|V_t\|^2
+C_\mu\delta^2b^{-1}g^4(1+G)^{-\frac32}.
\]
Applying the Gronwall-type inequality, we obtain the following estimate:
\[
\begin{aligned}
&GE_0
+2(1-\mu)\int_{t_0}^t bG\|V_t\|^2\,dr\\
&\le\Biggl[
\Biggl(G(t_0)E_0(t_0)
+C\int_{t_0}^t\left(gE_0
+\delta^2b^{-1}g^4
(1+G)^{-\frac32}\right)\,dr\Biggr)^{\frac12}
+\int_{t_0}^tG^{\frac12}\|F_2\|\,dr
\Biggr]^2.
\end{aligned}
\]
We apply Corollary~\ref{cor:E0} to the term
$\int_{t_0}^t gE_0\,dr$ on the right-hand side.
Moreover, the term involving $F_2$ can be estimated as
\[
\begin{aligned}
\left(\int_{t_0}^tG^{\frac12}\|F_2\|\,dr\right)^2
&\le C(\varepsilon+\delta\varepsilon)^2
\left(\int_{t_0}^t
g^{-1}(1+G)^{-\frac32}Q_\lambda^2\,dr\right)
\left(\int_{t_0}^t g\|V_x\|^2\,dr\right)\\
&\le C(\delta+\varepsilon)
\int_{t_0}^t g\|V_x\|^2\,dr.
\end{aligned}
\]
Thus, by taking $\varepsilon>0$ sufficiently small, we obtain the desired estimate.
\end{proof}

Differentiating \eqref{eq:V-equation} with respect to $x$, we obtain
\begin{equation}\label{eq:Vx}
V_{xtt}+\bigl(p'(V_x+\widetilde v)V_{xx}\bigr)_x+b(t)V_{xt}=F_{1x}
-\bigl(p'(V_x+\widetilde v)\widetilde v_x
-p'(\bar v)\bar v_x\bigr)_x=:F^{(1)}.
\end{equation}
We decompose $F^{(1)}$ as
$F^{(1)}=F_{1x}+F_2^{(1)}(V)+F_3^{(1)}$,
where
\[
\begin{aligned}
F_2^{(1)}(V)
&:=-\Bigl(\left(p'(V_x+\widetilde v)-p'(\widetilde v)\right)
\widetilde v_x\Bigr)_x \\
F_3^{(1)}
&:=-\Bigl(\left(p'(\widetilde v)-p'(\bar v)\right)\bar v_x
+p'(\widetilde v)\hat v_x\Bigr)_x.
\end{aligned}
\]
\eqref{eq:Vx}
also can be rewritten as
\begin{equation}\label{eq:Vx-g}
\bigl(g(t)V_{xt}+V_x\bigr)_t
+g(t)\bigl(p'(V_x+\widetilde v)V_{xx}\bigr)_x
=g(t)F^{(1)}.
\end{equation}

Differentiating \eqref{eq:V-equation} twice with respect to $x$, we obtain
\begin{equation}\label{eq:Vxx}
V_{xxtt}+\bigl(p'(V_x+\widetilde v)V_{xxx}\bigr)_x
+b(t)V_{xxt}
=F^{(2)},
\end{equation}
where
$F^{(2)}:=F_{1xx}+F_2^{(2)}(V)+\bigl(F_3^{(1)}\bigr)_x$ with
\[
F_2^{(2)}(V)
:=\bigl(F_2^{(1)}(V)\bigr)_x
-\Bigl(p''(V_x+\widetilde v)
\bigl(V_{xx}+\widetilde v_x\bigr)V_{xx}\Bigr)_x.
\]

\begin{definition}\label{def:higher-energies}
For $k=1,2$, we define the auxiliary energies and usual energies by
\[
\begin{aligned}
\widetilde E_k(t)
&:=-g^2\int p'(V_x+\widetilde v)
\bigl(\partial_x^{k+1}V\bigr)^2\,dx+\left\|g\partial_x^k\partial_tV
+\frac12\partial_x^kV\right\|^2
+\frac14\|\partial_x^kV\|^2,\\
E_k(t)
&:=-\int p'(V_x+\widetilde v)
\bigl(\partial_x^{k+1}V\bigr)^2\,dx
+\|\partial_x^k\partial_tV\|^2.
\end{aligned}
\]
\end{definition}

\begin{lemma}\label{lem:E1}
If $\varepsilon>0$ is sufficiently small, then the following estimate holds
on $[t_0,t_0+T]$:
\[
G\widetilde E_1
+\int_{t_0}^t gGE_1\,dr
\le C\left(
\|V(t_0)\|_{H^2}^2
+\|V_t(t_0)\|_{H^1}^2
+\delta^2Q_\lambda^2
\right).
\]
\end{lemma}

\begin{proof}
Differentiating $G(t)\widetilde E_1(t)$ and using
\eqref{eq:Vx-g}, we obtain
\[
\begin{aligned}
\frac{d}{dt}\bigl(G\widetilde E_1\bigr)
&\le g\widetilde E_1
-\left(\frac12-C\delta\right)gGE_1-g^2G\int p''(V_x+\widetilde v)
V_{xt}V_{xx}^2\,dx\\
&\quad +2gG\int F^{(1)}
\left(gV_{xt}+\frac12V_x\right)\,dx.
\end{aligned}
\]
Since
\begin{equation}\label{eq:E1-equivalence}
\widetilde E_1(t)
\asymp g(t)^2E_1(t)+\|V_x(t)\|^2,
\end{equation}
and
$g(t)^2/(1+G(t))\lesssim(1+t)^{-(1-\lambda)}$,
by taking $t_0$ sufficiently large, we have
\[
g(t)\widetilde E_1(t)
\le \frac14 g(t)G(t)E_1(t)+Cg(t)E_0(t)
\]
for all $t\ge t_0$.
Moreover, \eqref{eq:V-linfty} yields
\[
\begin{aligned}
g^2G\left|\int p''(V_x+\widetilde v)
V_{xt}V_{xx}^2\,dx\right|
&\le Cg^2G\|V_{xt}\|_{L^\infty}\|V_{xx}\|^2\\
&\le C\varepsilon \left(g^2Q_\lambda
(1+G)^{-\frac74}\right)gG\|V_{xx}\|^2\\
&\le C\varepsilon gGE_1.
\end{aligned}
\]

As in the proof of Lemma~\ref{lem:E0}, we have
\[
\begin{aligned}
g^2G\left|\int (F_{1x}+F^{(1)}_3)V_{xt}\,dx\right|
&\le \mu gG\|V_{xt}\|^2
+C_\mu\delta^2g^7(1+G)^{-\frac52}, \\
gG\left|\int (F_{1x}+F^{(1)}_3)V_x\,dx\right|
&\le \mu gG\|V_{xx}\|^2
+C_\mu\delta^2g^5(1+G)^{-\frac32}.
\end{aligned}
\]
For $F_2^{(1)}$, since
\[
|F_2^{(1)}|
\le C\left(
|\widetilde v_x||V_{xx}|
+\bigl(|\widetilde v_x|^2+|\widetilde v_{xx}|\bigr)|V_x|
\right),
\]
it follows that
\begin{equation}\label{eq:F2x-bound}
\|F_2^{(1)}\|
\le C\delta\left(
(1+G)^{-\frac12}\|V_{xx}\|
+(1+G)^{-1}\|V_x\|
\right).
\end{equation}
Thus, the term involving $F^{(1)}_2$ satisfies
\[
\begin{aligned}
gG\left|\int F_2^{(1)}
\left(gV_{xt}+\frac12V_x\right)\,dx\right|
&\le gG\|F^{(1)}_2\| \left( g\|V_{xt}\|+ \frac12 \|V_x\|\right)\\
&\le C\delta\bigl(g\|V_x\|^2+gGE_1\bigr).
\end{aligned}
\]
Therefore, we obtain
\[
\begin{aligned}
\frac{d}{dt}\bigl(G\widetilde E_1\bigr)
&\le-\left(\frac14-C(\delta+\varepsilon+\mu)\right)gGE_1
+CgE_0+C_\mu\delta^2g^5(1+G)^{-\frac32}.
\end{aligned}
\]
Integrating this inequality over $[t_0,t]$, we obtain the following
estimate:
\[
G\widetilde E_1
+\int_{t_0}^t gGE_1\,dr
\le C\Biggl(
G(t_0)\widetilde E_1(t_0)
+\delta^2Q_\lambda^2
+\int_{t_0}^t gE_0\,dr
\Biggr).
\]
Applying Corollary~\ref{cor:E0} to the third term on the right-hand
side, we obtain the desired estimate.
\end{proof}

\begin{lemma}\label{lem:E1-weighted}
If $\varepsilon>0$ is sufficiently small, then the following estimate holds
on $[t_0,t_0+T]$:
\[
G^2E_1
+\int_{t_0}^t bG^2\|V_{xt}\|^2\,dr
\le C\left(
\|V(t_0)\|_{H^2}^2
+\|V_t(t_0)\|_{H^1}^2
+\delta^2Q_\lambda^2
\right).
\]
\end{lemma}

\begin{proof}
Differentiating $G(t)^2E_1(t)$ and using
\eqref{eq:Vx} and
\eqref{eq:V-linfty}, we obtain
\[
\begin{aligned}
\frac{d}{dt}\bigl(G^2E_1\bigr)
&\le-2bG^2\|V_{xt}\|^2
+\bigl(2+C(\delta+\varepsilon)\bigr)gGE_1
+2G^2\int F^{(1)}(V)V_{xt}\,dx.
\end{aligned}
\]
The terms involving $F_{1x}$ and $F_3^{(1)}$ satisfy
\[
\begin{aligned}
G^2\left|\int\bigl(F_{1x}+F_3^{(1)}\bigr)V_{xt}\,dx\right|
&\le\mu bG^2\|V_{xt}\|^2
+C_\mu\delta^2b^{-1}g^4(1+G)^{-\frac32}.
\end{aligned}
\]
For the term involving $F_2^{(1)}$,
the estimate \eqref{eq:F2x-bound} gives
\[
\begin{aligned}
G^2\left|\int F_2^{(1)}V_{xt}\,dx\right|
&\le\mu bG^2\|V_{xt}\|^2
+C_\mu b^{-1}G^2\|F_2^{(1)}\|^2\\
&\le\mu bG^2\|V_{xt}\|^2
+C_\mu\delta\left(gG\|V_{xx}\|^2+g\|V_x\|^2\right).
\end{aligned}
\]
Here we used $b^{-1}\lesssim g$.
Therefore, integrating over
$[t_0,t]$ and applying Corollary~\ref{cor:E0} and Lemma~\ref{lem:E1}, we obtain
the desired estimate.
\end{proof}

The following two lemmas can be proved in the same way as
Lemmas~\ref{lem:E1} and \ref{lem:E1-weighted}.

\begin{lemma}\label{lem:E2}
If $\varepsilon>0$ is sufficiently small, then the following estimate holds
on $[t_0,t_0+T]$:
\[
G^2\widetilde E_2
+\int_{t_0}^t gG^2E_2\,dr
\le C\left(
\|V(t_0)\|_{H^3}^2
+\|V_t(t_0)\|_{H^2}^2
+\delta^2Q_\lambda^2
\right).
\]
\end{lemma}

\begin{proof}
Proceeding as in the proof of Lemma~\ref{lem:E1}, we obtain
\[
\begin{aligned}
\frac{d}{dt}\bigl(G^2\widetilde E_2\bigr)
&\le-\left(\frac14-C(\delta+\varepsilon+\mu)\right)gG^2E_2
+CgGE_1+C_\mu\delta^2g^5(1+G)^{-\frac32}\\
&\quad+2gG^2\int F_2^{(2)}
\left(gV_{xxt}+\frac12V_{xx}\right)\,dx.
\end{aligned}
\]
For $F_2^{(2)}$, we first note that
\[
\begin{aligned}
|F_2^{(2)}|
&\le C\Bigl[
\bigl(|\widetilde v_x|+|V_{xx}|\bigr)|V_{xxx}|+\bigl(
|\widetilde v_x|^2+|\widetilde v_{xx}|
+|\widetilde v_x||V_{xx}|+|V_{xx}|^2
\bigr)|V_{xx}|\\
&\quad+\bigl(
|\widetilde v_x|^3
+|\widetilde v_x||\widetilde v_{xx}|
+|\widetilde v_{xxx}|
\bigr)|V_x|
\Bigr].
\end{aligned}
\]
Hence, \eqref{eq:V-linfty} gives
\begin{equation}\label{eq:F2xx-bound}
\|F_2^{(2)}\|
\le C(\delta+\varepsilon)\Bigl[
(1+G)^{-\frac12}\|V_{xxx}\| +(1+G)^{-1}\|V_{xx}\|
+(1+G)^{-\frac32}\|V_x\|\Bigr].
\end{equation}
Then it follows that
\[
gG^2\left|\int F_2^{(2)}
\left(gV_{xxt}+\frac12V_{xx}\right)\,dx\right|
\le C(\delta+\varepsilon)
\bigl(gG^2E_2+gGE_1+gE_0\bigr).
\]
Integrating the resulting differential inequality over $[t_0,t]$ and
applying Corollary~\ref{cor:E0} and Lemma~\ref{lem:E1}, we obtain
the desired estimate.
\end{proof}

\begin{lemma}\label{lem:E2-weighted}
If $\varepsilon>0$ is sufficiently small, then the following estimate holds
on $[t_0,t_0+T]$:
\[
G^3E_2
+\int_{t_0}^t bG^3\|V_{xxt}\|^2\,dr
\le C\left(
\|V(t_0)\|_{H^3}^2
+\|V_t(t_0)\|_{H^2}^2
+\delta^2Q_\lambda^2
\right).
\]
\end{lemma}

\begin{proof}
Proceeding as in the proof of Lemma~\ref{lem:E1-weighted}, we obtain
\[
\begin{aligned}
\frac{d}{dt}\bigl(G^3E_2\bigr)
&\le-2(1-\mu)bG^3\|V_{xxt}\|^2
+\bigl(3+C(\delta+\varepsilon)\bigr)gG^2E_2\\
&\quad+C_\mu\delta^2b^{-1}g^4(1+G)^{-\frac32}
+2G^3\int F_2^{(2)}V_{xxt}\,dx.
\end{aligned}
\]
For the term involving $F_2^{(2)}$,
\eqref{eq:F2xx-bound} gives
\[
\begin{aligned}
G^3\left|\int F_2^{(2)}V_{xxt}\,dx\right|
&\le\mu bG^3\|V_{xxt}\|^2\\
&\quad+C_\mu(\delta+\varepsilon)\left(
gG^2\|V_{xxx}\|^2
+gG\|V_{xx}\|^2+g\|V_x\|^2
\right).
\end{aligned}
\]
Therefore, integrating over $[t_0,t]$ and applying Corollary~\ref{cor:E0} and
Lemmas~\ref{lem:E1}, \ref{lem:E2}, we obtain the desired estimate.
\end{proof}

If $\varepsilon>0$ is sufficiently small, the lemmas in this subsection
show that each quantity listed below is bounded by
$C\left(
\|V(t_0)\|_{H^3}^2
+\|V_t(t_0)\|_{H^2}^2
+\delta^2Q_\lambda(t)^2
\right)$:
\begin{center}
\renewcommand{\arraystretch}{2.0}
\begin{tabular}{c|l|c}
$(1)$
&$\displaystyle
\widetilde E_0(t)+\int_{t_0}^t g(r)E_0(r)\,dr$
&Corollary~\mbox{\ref{cor:E0}}\\[1ex] \hline
$(2)$
&$\displaystyle
G(t)E_0(t)+\int_{t_0}^t b(r)G(r)\|V_t(r)\|^2\,dr$
&Lemma~\mbox{\ref{lem:E0-weighted}}\\[1ex] \hline
$(3)$
&$\displaystyle
G(t)\widetilde E_1(t)+\int_{t_0}^t g(r)G(r)E_1(r)\,dr$
&Lemma~\mbox{\ref{lem:E1}}\\[1ex] \hline
$(4)$
&$\displaystyle
G(t)^2E_1(t)+\int_{t_0}^t b(r)G(r)^2\|V_{xt}(r)\|^2\,dr$
&Lemma~\mbox{\ref{lem:E1-weighted}}\\[1ex] \hline
$(5)$
&$\displaystyle
G(t)^2\widetilde E_2(t)+\int_{t_0}^t g(r)G(r)^2E_2(r)\,dr$
&Lemma~\mbox{\ref{lem:E2}}\\[1ex] \hline
$(6)$
&$\displaystyle
G(t)^3E_2(t)+\int_{t_0}^t b(r)G(r)^3\|V_{xxt}(r)\|^2\,dr$
&Lemma~\mbox{\ref{lem:E2-weighted}}\\[1ex]
\end{tabular}
\end{center}
In particular, estimates $(1)$, $(2)$, $(4)$, and $(6)$ imply that
\[
\sum_{k=0}^3G(t)^k\|\partial_x^kV(t)\|^2
\le C\left(
\|V(t_0)\|_{H^3}^2
+\|V_t(t_0)\|_{H^2}^2
+\delta^2Q_\lambda(t)^2
\right).
\]
The remaining estimates are derived in the next subsection, where we use
estimates $(2)$, $(4)$, and $(6)$ for
\begin{equation}\label{eq:dissipation}
\int_{t_0}^t b(r)G(r)^{k+1}
\|\partial_x^k\partial_tV(r)\|^2\,dr
\qquad k=0,1,2.
\end{equation}

\subsection{Energy estimates for $W$}\label{sec:energy-W}

In view of the estimates for \eqref{eq:dissipation}, we can aim to obtain the desired sharper estimates by deriving similar energy estimates for the time-differentiated equation. To this end, in this subsection, we derive energy estimates for $\partial_x^k W$ ($k=0,1$) from the equation for $W$ \eqref{eq:Wstar-equation}, which is formulated as the time-differentiated equation (see Subsection~\ref{sec:reformulations}). All estimates obtained in this subsection are summarized in a table at the end.

Rewriting \eqref{eq:Wstar-equation}, we obtain
\begin{equation}\label{eq:W}
W_{tt}+\bigl(p'(V_x^*+\widetilde v^*)W_x\bigr)_x
+\beta(t)W_t=R,
\end{equation}
where $\beta(t):=\frac{1+h'(t)}{h(t)}$ and $R:=R_1+R_2+R_3$ with
\[
\begin{aligned}
R_1
&:=\frac{1}{h(t)}\bigl(h(t)p(\bar v^*)_{xt}\bigr)_t,\\
R_2
&:=-\frac{1}{h(t)}\left(\bigl(p'(V_x^*+\widetilde v^*)
-p'(\widetilde v^*)\bigr)
\widetilde v_t^*\right)_x,\\
R_3
&:=-\frac{1}{h(t)}
\bigl(p(\widetilde v^*)-p(\bar v^*)\bigr)_{xt}.
\end{aligned}
\]
Here we note that
\[
\lim_{t\to\infty}\frac{\beta(t)}{b(t)}=1.
\]
Equation~\eqref{eq:W} can also be written as
\[
\bigl(h(t)W_t+W\bigr)_t
+h(t)\bigl(p'(V_x^*+\widetilde v^*)W_x\bigr)_x=h(t)R.
\]

\begin{definition}\label{def:J}
We define the auxiliary energy and usual energy by
\[
\begin{aligned}
\widetilde J_0(t)
&:=-h(t)^2\int p'(V_x^*+\widetilde v^*)W_x^2\,dx
+\left\|h(t)W_t+\frac12W\right\|^2
+\frac14\|W\|^2,\\
J_0(t)
&:=-\int p'(V_x^*+\widetilde v^*)W_x^2\,dx
+\|W_t\|^2.
\end{aligned}
\]
\end{definition}

\begin{lemma}\label{lem:J0}
If $\varepsilon>0$ is sufficiently small, then the following estimate holds
on $[t_0,t_0+T]$:
\[
\begin{aligned}
H^2\widetilde J_0
+\int_{t_0}^t hH^2J_0\,dr
&\le C\left(
\|V(t_0)\|_{H^3}^2+\|V_t(t_0)\|_{H^2}^2
+\delta^2Q_\lambda^2\right),\\
H^3J_0
+\int_{t_0}^t \beta H^3\|W_t\|^2\,dr
&\le C\left(
\|V(t_0)\|_{H^3}^2+\|V_t(t_0)\|_{H^2}^2
+\delta^2Q_\lambda^2\right).
\end{aligned}
\]
\end{lemma}

\begin{proof}
As in the proofs of Lemmas~\ref{lem:E1} and \ref{lem:E2}, we obtain
\[
\frac{d}{dt}\bigl(H^2\widetilde J_0\bigr)
\le-\left(\frac14-C(\delta+\varepsilon)\right)hH^2J_0
+ChH\|W\|^2+2hH^2\int R
\left(hW_t+\frac12W\right)\,dx,
\]
and, for the terms involving $R_1$ and $R_3$,
\[
\begin{aligned}
&hH^2\left|\int(R_1+R_3)
\left(hW_t+\frac12W\right)\,dx\right| \\
&\le \mu hH^2\bigl(\|W_t\|^2+\|W_x\|^2\bigr)
+C_\mu\delta^2h^5(1+H)^{-\frac32}.
\end{aligned}
\]
For $R_2$, since
\begin{equation}\label{eq:R2-bound}
\|R_2\|
\le C\delta\left(
(1+H)^{-1}\|V_{xx}^*\|
+(1+H)^{-\frac32}\|V_x^*\|
\right),
\end{equation}
we obtain
\[
\begin{aligned}
h^2H^2\left|\int R_2W_t\,dx\right|
&\le \mu hH^2\|W_t\|^2
+C_\mu\delta^2\left(
hH\|V_{xx}^*\|^2+h\|V_x^*\|^2
\right).
\end{aligned}
\]
Moreover, integration by parts gives
\[
hH^2\left|\int R_2W\,dx\right|
\le C\delta\left(
h\|V_x^*\|^2+hH^2\|W_x\|^2
\right).
\]
Integrating the resulting differential inequality over $[t_0,t]$ and
using $W=V_t^*/h$, estimates $(1)$--$(4)$ in
Subsection~\ref{sec:energy-V}, and Lemma~\ref{lem:tilde-difference}, we
obtain the first desired estimate.

By arguments analogous to those in the proofs of
Lemmas~\ref{lem:E1-weighted} and \ref{lem:E2-weighted}, together with
\eqref{eq:R2-bound}, we obtain the second desired estimate.
\end{proof}

Differentiating \eqref{eq:W} with respect to $x$, we obtain
\[
W_{xtt}+\bigl(p'(V_x^*+\widetilde v^*)W_{xx}\bigr)_x
+\beta(t)W_{xt}=R^{(1)},
\]
where $R^{(1)}:=R_{1x}+R_2^{(1)}+R_{3x}$ with
\[
R_2^{(1)}:=R_{2x}
-\Bigl(p''(V_x^*+\widetilde v^*)
\bigl(V_{xx}^*+\widetilde v_x^*\bigr)W_x\Bigr)_x.
\]

\begin{definition}\label{def:J1}
We define
\[
\begin{aligned}
\widetilde J_1(t)
&:=-h(t)^2\int p'(V_x^*+\widetilde v^*)W_{xx}^2\,dx
+\left\|h(t)W_{xt}+\frac12W_x\right\|^2
+\frac14\|W_x\|^2,\\
J_1(t)
&:=-\int p'(V_x^*+\widetilde v^*)W_{xx}^2\,dx
+\|W_{xt}\|^2.
\end{aligned}
\]
\end{definition}

\begin{lemma}\label{lem:J1}
If $\varepsilon>0$ is sufficiently small, then the following estimate holds
on $[t_0,t_0+T]$:
\[
\begin{aligned}
H^3\widetilde J_1
+\int_{t_0}^t hH^3J_1\,dr
&\le C\left(
\|V(t_0)\|_{H^3}^2+\|V_t(t_0)\|_{H^2}^2
+\delta^2Q_\lambda^2\right),\\
H^4J_1
+\int_{t_0}^t \beta H^4\|W_{xt}\|^2\,dr
&\le C\left(
\|V(t_0)\|_{H^3}^2+\|V_t(t_0)\|_{H^2}^2
+\delta^2Q_\lambda^2\right).
\end{aligned}
\]
\end{lemma}

\begin{proof}
We proceed as in the proof of Lemma~\ref{lem:J0}.
In doing so, we use
\[
\begin{aligned}
\|R_2^{(1)}\|
&\le C(\delta+\varepsilon)\Bigl(
(1+H)^{-\frac12}\|W_{xx}\|
+(1+H)^{-1}\|W_x\|\\
&\qquad+(1+H)^{-1}\|V_{xxx}^*\|
+(1+H)^{-\frac32}\|V_{xx}^*\|
+(1+H)^{-2}\|V_x^*\|
\Bigr).
\end{aligned}
\]
This estimate, together with estimates $(1)$, $(3)$, and $(5)$ in
Subsection~\ref{sec:energy-V} and Lemma~\ref{lem:J0}, yields the desired
estimates.
\end{proof}

If $\varepsilon>0$ is sufficiently small, the lemmas in this subsection
show that each quantity listed below is bounded by
$C\left(
\|V(t_0)\|_{H^3}^2
+\|V_t(t_0)\|_{H^2}^2
+\delta^2Q_\lambda(t)^2
\right)$:
\begin{center}
\renewcommand{\arraystretch}{2.0}
\begin{tabular}{c|l|c}
$(7)$
&$\displaystyle
H(t)^2\widetilde J_0(t)
+\int_{t_0}^t h(r)H(r)^2J_0(r)\,dr$
&Lemma~\mbox{\ref{lem:J0}}\\[1ex] \hline
$(8)$
&$\displaystyle
H(t)^3J_0(t)
+\int_{t_0}^t \beta(r)H(r)^3\|W_t(r)\|^2\,dr$
&Lemma~\mbox{\ref{lem:J0}}\\[1ex] \hline
$(9)$
&$\displaystyle
H(t)^3\widetilde J_1(t)
+\int_{t_0}^t h(r)H(r)^3J_1(r)\,dr$
&Lemma~\mbox{\ref{lem:J1}}\\[1ex] \hline
$(10)$
&$\displaystyle
H(t)^4J_1(t)
+\int_{t_0}^t \beta(r)H(r)^4\|W_{xt}(r)\|^2\,dr$
&Lemma~\mbox{\ref{lem:J1}}\\[1ex]
\end{tabular}
\end{center}

In particular, estimates $(7)$, $(8)$, and $(10)$ imply that
\[
\sum_{k=0}^2 g(t)^{-2}G(t)^{k+2}
\|\partial_x^k\partial_tV(t)\|^2
\le C\left(
\|V(t_0)\|_{H^3}^2
+\|V_t(t_0)\|_{H^2}^2
+\delta^2Q_\lambda(t)^2
\right).
\]
The estimates obtained in Subsections~\ref{sec:energy-V} and
\ref{sec:energy-W} complete the proof of Proposition~\ref{prop:apriori}.

\subsection{Finite-time estimates and conclusion}\label{sec:conclusion}

The preceding arguments establish a priori estimates for solutions of the reformulated problem~\eqref{eq:V-equation} on $[t_0,t_0+T]$.
In this subsection, for the fixed $t_0$, we construct a small
solution on $[0,t_0]$.  For this purpose, we prepare the following lemma.

\begin{lemma}\label{lem:finite-time-energy}
There exist constants $\varepsilon_0>0$ and $C=C(t_0)>0$ such that if $\delta\le\varepsilon_0$ holds and
$V\in\bigcap_{i=0}^2 C^i\bigl([0,T];H^{3-i}(\mathbb{R})\bigr)$
$(0<T\le t_0)$
is a solution of \eqref{eq:V-equation} satisfying
\[
\sup_{0\le t\le T}
\left(\|V(t)\|_{H^3}^2+\|V_t(t)\|_{H^2}^2\right)
\le\varepsilon_0^2,
\]
then the following holds:
\[
\sup_{0\le t\le T}
\left(\|V(t)\|_{H^3}^2+\|V_t(t)\|_{H^2}^2\right)
\le C\left(\|V_0\|_{H^3}^2+\|V_1\|_{H^2}^2+\delta^2\right).
\]
\end{lemma}

\begin{proof}
Differentiating $\widetilde E_0(t)$, $E_0(t)$, $E_1(t)$, and $E_2(t)$ as in
Subsection~\ref{sec:energy-V} and applying the standard Gronwall inequality,
we successively obtain
\[
\begin{aligned}
\widetilde E_0(t)
&\le e^{Ct_0}\left(
\widetilde E_0(0)+\int_0^tE_0(r)\,dr+\delta^2\right),\\
E_0(t)
&\le e^{Ct_0}\bigl(E_0(0)+\delta^2\bigr),\\
E_1(t)
&\le e^{Ct_0}\bigl(E_0(0)+E_1(0)+\delta^2\bigr),\\
E_2(t)
&\le e^{Ct_0}
\bigl(E_0(0)+E_1(0)+E_2(0)+\delta^2\bigr),
\end{aligned}
\]
for $0\le t\le T$.
These estimates give the desired result.
\end{proof}

\begin{remark}
By Lemma~\ref{lem:finite-time-energy}, if
$\|V_0\|_{H^3}+\|V_1\|_{H^2}+\delta$ is sufficiently small, a solution can
be constructed on $[0,t_0]$ while keeping
$\|V(t)\|_{H^3}^2+\|V_t(t)\|_{H^2}^2$ small throughout this interval.
\end{remark}

Combining Lemma~\ref{lem:finite-time-energy} with the local existence, we construct a sufficiently small solution on $[0,t_0]$. 
Then combining the a priori estimate (Proposition~\ref{prop:apriori}) with local existence and applying a standard continuation argument, we obtain Theorem~\ref{thm:global-existence} and hence Theorem~\ref{thm:main}.

\makeatletter
\renewcommand{\@biblabel}[1]{\textbf{[#1]}}
\makeatother

\end{document}